\documentclass[10pt]{amsart}
\usepackage{amsmath,amscd}
\usepackage{amsbsy}
\usepackage{amssymb}
\usepackage{amscd,amsthm}
\usepackage[all,cmtip]{xy}
\usepackage[russian,english]{babel}
\usepackage[utf8]{inputenc}
\usepackage{xcolor}

\usepackage{hyperref}

\renewcommand{\epsilon}{\varepsilon}

\renewcommand{\ell}{x}

\newtheorem{thm}{Theorem}\numberwithin{thm}{section}
\newtheorem{lem}[thm]{Lemma}
\newtheorem{prop}[thm]{Proposition}
\newtheorem{cor}[thm]{Corollary}

\newtheorem{rema}[thm]{Remark}

\newtheorem*{con2}{Conjecture}

\begin{document}
	\begin{center}
		\huge{On a generalization of the Brocard--Ramanujan Diophantine equation}\\[1cm] 
	\end{center}
	\begin{center}
	\large{Sa$\mathrm{\check{s}}$a Novakovi$\mathrm{\acute{c}}$}\\[0,5cm]
		{\small February 2026}\\[0,2cm]
	\end{center}
	\begin{center}
		\emph{dedicated to my father.}\\[0,3cm]
	\end{center}
	{\small \textbf{Abstract}. 
		Let $Q_1,...,Q_r\in \mathbb{Z}[x]$ be polynomials having $0$ as a root. Let $f(x,y)\in\mathbb{Z}[x,y]$ be a homogeneous polynomial with factorization $f(x,y)=f_1(x,y)^{e_1}\cdots f_u(x,y)^{e_u}$,
		where $f_i(x,y)$ are irreducible homogeneous polynomials of degree $d_i\geq 2$. Fix some positive integers $A_1,...,A_r$. We show that under certain conditions, the diophantine equation $\prod_{i=1}^rQ_i(A_i^{n_i}n_i!)=f(x,y)$ has finitely many integer solutions.}
	\begin{center}
		\tableofcontents
	\end{center}
	\section{Introduction}
	The theory of Diophantine equations has a long and rich history and has attracted the attention of many mathematicians. In particular, the study of diophantine equations involving factorials have been studied extensively. For example Brocard \cite{BR}, and independently Ramanujan \cite{RA}, asked to find all integer solutions for $n!=x^2-1$. It is still an open problem, known as Brocard's problem, and it is believed that the equation has only three solutions $(x,n)=(5,4), (11,5)$ and $(71,7)$. Calculations up to $n=10^9$ confirmed this \cite{BG}. Overholt \cite{O} observed that a weak form of Szpiro's conjecture implies that the Brocard--Ramanujan equation has finitely many integer solutions. Quite recently, Naciri \cite{NA} showed that  $n!+1=x^2$ has only finitely many integer solutions, assuming $x\pm1$ is a $k$-free integer or a prime power. Some further examples of similar equations are:
	\begin{itemize}
		\item[1)] $n!+A=x^2$, see \cite{DAB} and \cite{DK},
		\item[2)] $n!=p(x)$, where $p(x)\in\mathbb{Z}[x]$, see \cite{L},
		\item[3)] $\overset{n}{\underset{\mathrm{gcd}(k,n)\neq k}{\prod}}=p(x)$, where $p(x)\in \mathbb{Z}[x]$, see \cite{KL},
		\item[4)] $n!=f(x,y)$, where $f(x,y)\in \mathbb{Z}[x,y]$ is a homogeneous polynomial, see \cite{WT}
		\item[5)] $Q(n!)=\overset{l}{\underset{i=1}{\prod}}p_i(x_i,y_i)^{\alpha_i}$, where $p_i$ are homogenous polynomials of degree $\geq 2$ and $Q$ a polynomial having $0$ as a root with odd multiplicity, see \cite{NA1}.
	\end{itemize}
	For equation 1), it was shown that the number of integer solutions is finite, provided weak form of Szpiro's or the Hall conjecture holds and $A$ is a square. The equation in 2) has finitely many integer solutions, provided the abc conjecture holds. The same was proven for 3). Takeda \cite{WT} studied 4) and proved that if $f$ is irreducible and homogeneous of degree $\geq 2$, there are finitely many factorials that can be represented by $f(x,y)$. And recently, Naciri \cite{NA1} showed that 5) has only finitely many integer solutions $(n,x_1,...,x_l,y_1,...,y_l,\alpha_1,...,\alpha_l)$. There are a lot of more diophantine equations involving factorials and polynomials that have been studied and we refer the interested reader to \cite{BN}, \cite{NO} and the references therein. There are also variants of the Brocard--Ramanujan equation that have been studied intensively in the literature. We do not want to give a complete list here and refer for instance to \cite{KL} and \cite{NA1}. At this point, we also want to mention the paper of Berend and Harmse \cite{BH1}, where the authors studied diophantine equations of the form $P(x)=H_n$, where $H_n$ are certain divisible sequences, namely $H_n=n!$ or $H_n=p_n\#$ or $H_n=[1,...,n]$. In \cite{NO} and \cite{NOO} the author considered the cases where $H_n=bA^nn!$ and $H_n=bA^nn!!$ and showed that under the abc conjecture the equations $P(x)=\prod_{i=1}^rA_i^{n_i}n_i!$ and $P(x)=\prod_{i=1}^rA_i^{n_i}n_i!!$ have only finitely many integer solutions $(x,n_1,...,n_r)$. Furthermore, the author also studied, more generally, the equation $P(x)=\prod_{i=1}^rH_{n_i}$, where $H_{n_i}$ are arbitrary divisible sequences (see \cite{NO2}). In the present note, we focus on the equation 5), treated by Naciri \cite{NA1}. In fact, we want to study slight modifications. Our main results are:
	
	\begin{thm}
			Let $f(x,y)=a_dx^d+a_{d-1}x^{d-1}y+\cdots + a_{1}xy^{d-1}+a_0y^d\in \mathbb{Z}[x,y]$ be an irreducible homogeneous polynomial and let $Q_1,...,Q_r\in\mathbb{Z}[x]$ be polynomials having $0$ as a root of multiplicity $l_i>0$. Fix positive integers $A_1,...,A_r$ and assume $d>l_1+\cdots +l_r\geq 1$. Then there are finitely many $(n_1,...,n_r)$ such that $\prod_{i=1}^rQ_i(A_i^{n_i}n_i!)$ is represented by $f(x,y)$. If $l_1+\cdots +l_r\geq 2$, then the  diophantine equation
			$$
			\prod_{i=1}^rQ_i(A_i^{n_i}n_i!)=f(x,y)
			$$
			has only finitely many integer solutions $(n_1,...,n_r,x,y)$.
	\end{thm}
	\begin{cor}
			Let $Q_i$ and $A_i$ be as in Theorem 1.1. For any irreducible $f(x)\in\mathbb{Z}[x]$ with degree $d>l_1+\cdots+l_r$, the equation 
				$$
			\prod_{i=1}^rQ_i(A_i^{n_i}n_i!)=f(x)
			$$
			 has only finitely many integer solutions $(n_1,...,n_r,x)$.
	\end{cor}
	\begin{thm}
			Let $f(x,y)\in\mathbb{Z}[x,y]$ be a homogeneous polynomial with factorization 
		\begin{center}
			$f(x,y)=f_1(x,y)^{e_1}\cdots f_u(x,y)^{e_u}$,
		\end{center}
		where $f_i(x,y)$ are irreducible homogeneous polynomials of degree $d_i\geq 2$. Furthermore, let $Q_1,...,Q_r\in\mathbb{Z}[x]$ be polynomials having $0$ as a root of multiplicity $l_i>0$. Fix some positive integers $A_1,...,A_r$ and assume $d_1e_1+\cdots +d_ue_u>l_1+\cdots +l_r$ or $\mathrm{min}\{d_1e_1,...,d_ue_u\}>l_1+\cdots+l_r$,
		Then there are finitely many $(n_1,...,n_r)$ such that $\prod_{i=1}^rQ_i(A_i^{n_i}n_i!)$ is represented by $f(x,y)$. If $l_1+\cdots +l_r\geq 2$, then the diophantine equation
		$$
		\prod_{i=1}^rQ_i(A_i^{n_i}n_i!)=f(x,y)
		$$
		has only finitely many integer solutions $(n_1,...,n_r,x,y)$.
	\end{thm}
	\noindent
	We now explain the relation to \cite{NA1}, Theorem 1. In fact, if all $d_i$ are assumed to be even and if $r=1$ and $l_1$ is assumed to be odd, the proof of Theorem 1.3 still works. And since the $f_i$ do not change signs, we even get:
\begin{thm}
Let $d_i$ be even numbers and $f_i(x,y)$ irreducible homogeneous polynomials of degree $d_i$. Furthermore, let $Q_1\in\mathbb{Z}[x]$ be a polynomial having $0$ as a root of odd multiplicity $l_1$. Fix some positive integer $A_1$. Then the diophantine equation
$$
Q_1(A_1^{n_1}n_1!)=f_1(x,y)^{e_1}\cdots f_u(x,y)^{e_u}
$$
has only finitely many integer solutions $(n_1,x,y,e_1,...,e_u)$.
\end{thm}
\noindent
In the context of Corollray 1.2, we can also consider the case where $f(x)$ is reducible. We can use the abc conjecture to show that $\prod_{i=1}^rQ_i(A_i^{n_i}n_i!)=f(x)$ has finitely many integers solutions. We recall the abc conjecture.
For a non-zero integer $a$, let $N(a)$ be the \emph{algebraic radical}, namely $N(a)=\prod_{p|a}{p}$. Note that 
\begin{eqnarray*}
	N(a)=\prod_{p|a}{p}\leq \prod_{p\leq a}{p}< 4^a,
\end{eqnarray*}
where the last inequality follows from a Chebyshev-type result in elementary prime number theory and is called the Finsler inequality. 
\begin{con2}[abc conjecture]
	For any $\epsilon >0$ there is a constant $K(\epsilon)$ depending only on $\epsilon$ such that whenever $a,b$ and $c$  are three coprime and non-zero integers with $a+b=c$, then 
	\begin{eqnarray}
		\mathrm{max}\{|a|,|b|,|c|\}<K(\epsilon)N(abc)^{1+\epsilon}
	\end{eqnarray}
	holds.
\end{con2} 

\begin{thm}
	Let $f(x)\in\mathbb{Q}[x]$ be a polynomial of degree $d\geq 2$ which is not monomial and has at least two distinct roots and let $Q_i$ and $A_i$ be as in Theorem 1.1. Furthermore, assume $f(x)\neq \prod_{i=1}^rQ_i(x)$. Then the abc conjecture implies that the equation
	$$
	\prod_{i=1}^rQ_i(A_i^{n_i}n_i!)=f(x)
	$$
	has only finitely many integers solutions $(n_1,...,n_r,x)$.
\end{thm}
\begin{rema}
	\textnormal{In the case $f(x)=Q(x)$, one obviously has infinitely many integers solutions for $f(x)=Q(A^nn!)$ simply by setting $x=A^nn!$. If $f(x)= \prod_{i=1}^rQ_i(x)$ and if, for instance, $A_1=A_2=\cdots =A_r$, one can set $n_1=n_2=\cdots =n_r$ and $x=A_1^{n_1}n_1!$ to find infinitely any solutions for $\prod_{i=1}^rQ_i(A_i^{n_i}n_i!)=f(x)$.} 
\end{rema}
\begin{prop}
	Let $Q_i$ and $A_i$ be as in Theorem 1.1 and assume that $d>l_1+\cdots +l_r$. Then the equation $x^d=\prod_{i=1}^rQ_i(A_i^{n_i}n_i!)$ has finitely many integer solutions $(n_1,...,n_r,x)$.
\end{prop}
	\section{Preliminary Results}
	\noindent
We follow the notation of \cite{WT}. Let $f(x,y)=a_dx^d+a_{d-1}x^{d-1}y+\cdots + a_{d-1}xy^{d-1}+a_0y^d$ be an irreducible polynomial and let $K_F$ be the splitting field of $f(x,1)$. Let $G_F=\mathrm{Gal}(K_F/\mathbb{Q})$ and let $p$ be a prime and $\textbf{p}$ a prime ideal of $\mathcal{O}_{K_F}$ lying above $p$. Consider the Galois group $\mathrm{Gal}((\mathcal{O}_{K_F}/\textbf{p})/(\mathbb{Z}/p\mathbb{Z}))$. It is known that this group is cyclic with unique generator $\sigma: x\mapsto x^p$. Then the Frobenius map $(p,K_F/\mathbb{Q})$ of $p$ is in the image of $\sigma$ in the Galois group $G_F$. If $(p,K_F/\mathbb{Q})$ belongs to a conjugacy class $C$ of $G_F$, we say $p$ corresponds to $C$ (see \cite{WT}, chapter 2 for details). We denote the set of primes corresponding to some $C$ by $P(C)$. Denote by $C_F$ the set of conjugacy classes of the Galois group $G_F=\mathrm{Gal}(K_F/\mathbb{Q})$ whose cycle type $[h_1,...,h_s]$ satisfies $h_i\geq 2$. For a cycle $\sigma$, the cycle type is defined as the ascending ordered list $[h_1,...,h_s]$ of the sizes of cycles in the cycle decomposition of $\sigma$. For further details we refer to Chapters 2, 3 and 4 in \cite{WT}. 
Let $g=\mathrm{gcd}(a_d,...,a_0)$ and $N=gp_1\cdots p_uq_1^{l_1}\cdots q_v^{l_v}$, where $q_i$ are primes corresponding to a conjugacy class in $C_F$ satisfying $\mathrm{gcd}(q_i,a\Delta_{mod})=1$ with $a\in\{a_d,a_0\}$ and $p_j$ are the other primes (see \cite{WT}, Lemma 3.1 for details).  Here $\Delta_{mod}$ denotes a certain modified discriminant of $f(x,1)$ which is defined by 
\begin{center}
	$\Delta_{mod}=\frac{\Delta_{f(x,1)}}{\mathrm{gcd}(a_d,...,a_0)^{2d-2}}$.
\end{center}
\begin{lem}[\cite{WT}, Lemma 3.1]
		Let $f(x,y)=a_dx^d+a_{d-1}x^{d-1}y+\cdots + a_{1}xy^{d-1}+a_0y^d\in \mathbb{Z}[x,y]$ be an irreducible homogeneous polynomial and $g=\mathrm{gcd}(a_d,...,a_0)$. Let $N$ be an integer with $N=gp_1\cdots p_sq_1^{v_1}\cdots q_t^{v_t}$, where $q_i$ are primes corresponding to a conjugacy class $C\in C_F$ satisfying $\mathrm{gcd}(q_i,a\Delta_{mod})=1$ and $p_i$ are other primes. If $N$ is represented by $f(x,y)$, then $d\mid v_j$ for $j=1,...,t$.  
\end{lem}
\begin{lem}[\cite{WT}, Lemma 3.2]
		Let $f(x,y)\in\mathbb{Z}[x,y]$ be a homogeneous polynomial of degree $d$ with irreducible factorization 
	$$
		f(x,y)=\prod_{i=1}^uf_i(x,y),
		$$
and $g=\mathrm{gcd}(a_d,...,a_0)$. Let $N$ be an integer with $N=gp_1\cdots p_sq_1^{v_1}\cdots q_t^{v_t}$, where $q_i$ are primes corresponding to a conjugacy class $C\in \cap_{j=1}^uC_{F_j}$ with $\mathrm{gcd}(q_i,a\Delta_{mod})=1$ and $p_i$ are other primes. If $N$ is represented by $f(x,y)$, then $d\mid v_j$ for $j=1,...,t$.  
\end{lem}
We consider prime ideals $\textbf{q}$ corresponding to a conjucacy class $C$ such that its ideal norm $\mathfrak{R}\textbf{q}$ is of the form $q^f$. 
\begin{thm}[\cite{WT}, Theorem 3.6]
	Let $L$ be the Galois closure of $K/\mathbb{Q}$ with $[L:\mathbb{Q}]=k$ and $p$ a prime corresponding to a conjugacy class $C$ of $\mathrm{Gal}(L/\mathbb{Q})$. For any $A>1$ there is an effective computable constant $c(A)>0$ such that for $p^{f}>\mathrm{exp}(c(A)k(\mathrm{log}(D_L))^2)$ there exists a prime ideal $\textbf{q}$ with ideal norm $\mathfrak{R}\textbf{q}=q^f\in (p^{f},Ap^{f})$, where $q$ corresponds to $C$. 
\end{thm}
\begin{lem}[\cite{WT}, Lemma 2.1]
	Let $H$ be a trasitive subgroup of $S_n$ of degree $n\geq 2$. Then there exists an element $\sigma\in H$ such that $\sigma(i)\neq i$ for all $i=1,...,n$.
\end{lem}
\noindent
Finally, we need the following result.
\begin{thm}[\cite{NO2}, Theorem 2.3]
	Let $P(x)\in\mathbb{Q}[x]$ be a polynomial of degree $d\geq 2$ which is not monomial and has at least two distinct roots. Let $F(n_1,...,n_r)$ be a function and assume there exists an $\epsilon>0$, such that $N(F(n_1,...n_r))^{1+\epsilon}=o(F(n_1,...n_r))$ as $(n_1,...,n_r)\rightarrow \infty$. Then the abc conjecture implies that $P(x)=F(n_1,...,n_r)$ has finitely many integer solutions $(n_1,...,n_r,x)$.
\end{thm}
\begin{rema}
	\textnormal{We want to mention that in Theorem 2.5 it is implicitely assumed that $r\geq 2$. Otherwise, one could set $F(n)=P(n)$ and would in this way obtain infinitely many solutions $x=n$. An argument can be read of from the proof of \cite{L}, Proposition 1, as it is the same argument that can be used to prove Theorem 2.5.}
\end{rema}
	
	\section{Proof of Theorems 1.1 and 1.3}
	\noindent
	We first prove Theorem 1.1 Since $Q_i$ are polynomials having $0$ as a root of multiplicity $l_i>0$, we can write
	$$
	Q_i(x)=x^{l_i}(f_{i0}+f_{i1}x+\cdots f_{ic_i}x^{c_i}),
	$$
	with $c_i\in \mathbb{N}$, $f_{i0},...,f_{ic_i}\in \mathbb{Z}$ and $f_{i0}\neq 0$. We set $f_{i0}+f_{i1}x+\cdots f_i{c_i}x^{c_i}:= R_i(x)$ Then
	$$
	\prod_{i=1}^rQ_i(A_i^{n_i}n_i!)=\prod_{i=1}^r(A_i^{n_i}n_i!)^{l_i}\cdot \prod_{i=1}^rR_i(A_i^{n_i}n_i!).
	$$
	Our diophantine equation then becomes
	$$
	\prod_{i=1}^r(A_i^{n_i}n_i!)^{l_i}\cdot \prod_{i=1}^rR_i(A_i^{n_i}n_i!)=f(x,y).
	$$
	The assumption $d\geq 2$ and Lemma 2.1 then imply that if $N$ is represented by $f(x,y)$ and $q|N$ for a prime $q$ corresponding to $C\in C_F$ satisfying $\mathrm{gcd}(q,a\Delta_{mod})=1$, then $N$ is divisible by $q^d$ at least. Notice that Lemma 2.4 guaranties that $C_F\neq \emptyset$. Without loss of generality, let $n_r\geq n_{r-1}\geq \cdots \geq n_1$ and assume $q$ is big enough, say $q>8\mathrm{max}\{A_1,...,A_r, f_{ij}\}$. Such primes exist because the set of such primes has positive density according to Chebotarjev density theorem. Since the second smallest positive integer divisible by $q$ is $2q$, there is no solution for our equation $\prod_{i=1}^u(A_i^{n_i}n_i!)^{l_i}\cdot \prod_{i=1}^uR_i(A_i^{n_i}n_i!)=f(x,y)$ if $q<n_r<2q$, because
	$$
	\nu_q(\prod_{i=1}^r(A_i^{n_i}n_i!)^{l_i}\cdot \prod_{i=1}^rR_i(A_i^{n_i}n_i!))\leq \sum_{i=1}^rl_i<d.
	$$
	Notice that $d>l_1+\cdots +l_r$ by assumption and that 
	$$
	\nu_p((A_i^{n_i}n_i!)^{l_i}\cdot R_i(A_i^{n_i}n_i!))\leq l_i, 
	$$ 
	since $p$ divides neither $A_i$ nor $f_{i0}$. Now apply Theorem 2.3 (as in the proof of Theorem 4.1 of \cite{WT}). Let $\alpha$ be a root of $f(x,1)$ and let $k$ be the degree of $K_F/\mathbb{Q}$. We denote the ring of integers by $\mathcal{O}_{\alpha}$. Let $\textbf{q}$ be a prime ideal of $\mathcal{O}_{\alpha}$ corresponding to $C$ with $\mathfrak{R}\textbf{q}= q^f>\mathrm{max}(\mathrm{exp}(c(A)k(\mathrm{log}(D_L))^2), (a_d\Delta_{mod})^f,(a_0\Delta_{mod})^f)$. Theorem 2.3 shows that there exists a $\textbf{q}'$ corresponding to $C$ with $\mathfrak{R}\textbf{q}'=q'^f\in (q^f,2q^f)$. Hence there is a prime $q'$ corresponding to $C$ such that $q'\in(q,2q)$ satisfying $\mathrm{gcd}(q',a\Delta_{mod})=1$. Therefore, by the same argument as before, if $q'<n_r<2q'$ there are no integer solutions for  $\prod_{i=1}^r(A_i^{n_i}n_i!)^{l_i}\cdot \prod_{i=1}^rR_i(A_i^{n_i}n_i!)=f(x,y)$. By an induction argument we conclude that whenever $n_r>q$ there are no integer solutions. This shows that there are only finitely many $(n_1,...,n_r)$ such that $\prod_{i=1}^u(A_i^{n_i}n_i!)^{l_i}\cdot \prod_{i=1}^uR_i(A_i^{n_i}n_i!)$  is represented by $f(x,y)$. Now if $l_1+\cdots +l_r\geq 2$, we have $d\geq 3$ and we can use Thue's theorem to conclude that there are indeed only finitely many integer solutions.\\
	\noindent
	We now prove Theorem 1.3. If all $d_i\geq 2$, then proceed as follows: let $K_{F_j}$ be the splitting field of $f_j(x,1)$ and denote by $C_{F_j}$ be the set of conjugacy classes whose cycle type has sizes $\geq 2$. Now let $q$ be a prime corresponding to a conjugacy class $C\in \cap_{j=1}^uC_{F_j}$ satisfying $\mathrm{gcd}(q,a\Delta_{mod})=1$. Assume $N$ is represented by $f(x,y)$ and $q\mid N$. Therefore, there are $x_0,y_0$ such that $q|f(x_0,y_0)$. Then there exists a polynomial $f_j(x,y)=a_{j,d_j}x^{d_j}+\cdots + a_{j,0}y^{d_j}$ such that $q|f_j(x_0,y_0)$. Now $\mathrm{gcd}(q,a\Delta_{mod})=1$ with  $a\in\{a_d,a_0\}$ implies $\mathrm{gcd}(q,a(j)\Delta_{j, mod})=1$ where $a(j)\in\{a_{j,d_j},a_{j,0}\}$. Here $\Delta_{j, mod}$ denotes the modified discriminant of $f_j(x,1)$. It follows from Lemma 2.1 that $q^d|f(x,y)$. As in the proof of Theorem 1.1 above, we may assume $n_r\geq \cdots \geq n_1$ and can choose $q$ big enough. Again, one concludes the statement from an inductions argument. The rest of the proof is the same as in the proof of Theorem 1.1. Now let $\mathrm{min}\{d_1e_1,...,d_ue_u\}=d_{i_0}e_{i_0}$ and assume $d_{i_0}=1$ and $e_{i_0}>1$. Then one can argue as in \cite{WT}, Theorem 4.7 to conclude that, let's say, $q$ divides $f(x,y)$ at least $e_{i_0}>1$ times. Again, the rest of the proof is as in the proof of Theorem 1.1 above.
	\begin{rema}
		\textnormal{For the proof of Theorem 1.3 one can use Lemma 2.2 directly. In the proof above, we made the argument explicit. In fact, we gave an argument why Lemma 2.2 follows form Lemma 2.1.}
		\end{rema}
	\section{Proof of Theorem 1.5 and Proposition 1.7}
	\noindent
	To keep notation as simple as possible, we prove the statement for $r=2$, because the general case is completely analogous. So we consider the equation
	$$
	Q_1(A^nn!)\cdot Q_2(B^mm!)=f(x).
	$$
	Let $F(n,m)=Q_1(A^nn!)\cdot Q_2(B^mm!)$. According to Theorem 2.5, we have to show that
	$$
	\frac{N(F(n,m))^{1+\epsilon}}{F(n,m)}\rightarrow 0
	$$
	for a suitable $\epsilon >0$ as $(n,m)\rightarrow \infty$. Denote the multiplicity of zero for $Q_1$ by $p$ and that for $Q_2$ by $q$. In the proof of Theorem 1.1, we obtained
	$$
	Q_1(A^nn!)=(A^nn!)^pR_1(A^nn!) \quad \textnormal{and} \quad Q_2(B^mm!)=(B^mm!)^qR_2(B^mm!),
	$$
	where $R_1,R_2$ are integer polynomials of the form
	$$
	R_1(x)=a_0+a_1x+\cdots +a_ex^e \quad \textnormal{and} \quad R_2(x)=b_0+b_1x+\cdots +b_fx^f,
	$$
	with $a_0\neq 0$ and $b_0\neq 0$. We have
	$$
	N(F(n,m))\leq N((A^nn!)^p)N((B^mm!)^q)N(R_1(A^nn!))N(R_2(B^mm!)).
	$$
	Since $N((A^nn!)^p)\leq N(A)4^n$ and $N((B^mm!)^q)\leq N(B)4^m$, we get
	$$
	N(F(n,m))\leq C4^{n+m}\cdot N(R_1(A^nn!))N(R_2(B^mm!)), 
	$$
	where $C=N(A)N(B)$. Since $N(x)\leq x$ for any positive integer $x$, we moreover have 
	$$
	N(F(n,m))\leq C4^{n+m}\cdot R_1(A^nn!)\cdot R_2(B^mm!).
	$$
	This yields
	$$
	N(F(n,m))^{1+\epsilon}\leq C'4^{(n+m)(1+\epsilon)}\cdot R_1(A^nn!)^{1+\epsilon}\cdot R_2(B^mm!)^{1+\epsilon}
	$$
	and therefore
	$$
	\frac{N(F(n,m))^{1+\epsilon}}{F(n,m)}\leq \frac{C'4^{(n+m)(1+\epsilon)}}{(A^nn!)^p(B^mm!)^q}\frac{R_1(A^nn!)^{1+\epsilon}\cdot R_2(B^mm!)^{1+\epsilon}}{R_1(A^nn!)\cdot R_2(B^mm!)}.
	$$
	The right hand side can be simplyfied and we obtain
	$$
	\frac{N(F(n,m))^{1+\epsilon}}{F(n,m)}\leq \frac{C'4^{(n+m)(1+\epsilon)}}{(A^nn!)^p(B^mm!)^q}\cdot R_1(A^nn!)^{\epsilon}\cdot R_2(B^mm!)^{\epsilon}.
	$$
	For a polynomial $g(x)\in \mathbb{Z}[x]$ there is allways a constant $D>0$ such that $g(l)\leq D\cdot l^{\mathrm{deg}(g)}$ for all positive integers $l$. The constant $D$ can be chosen to be the sum of the absolute values of all coefficients of $g$. We apply this to the polynomials $R_1$ and $R_2$ and get
	$$
	\frac{N(F(n,m))^{1+\epsilon}}{F(n,m)}\leq \frac{C'4^{(n+m)(1+\epsilon)}}{(A^nn!)^p(B^mm!)^q}\cdot D_1^{\epsilon}(A^nn!)^{\mathrm{deg}(R_1)\epsilon}\cdot D_2^{\epsilon}(B^mm!)^{\mathrm{deg}(R_2)\epsilon}.
	$$
	If we set $C''=C'(D_1D_2)^{\epsilon}$ and $d_i=\mathrm{deg}(R_i)$, we can rewrite the right hand side of the inequality and find
	$$
	\frac{N(F(n,m))^{1+\epsilon}}{F(n,m)}\leq C''4^{(n+m)(1+\epsilon)}\frac{ A^{n\cdot d_1\cdot \epsilon}\cdot B^{m\cdot d_2\cdot \epsilon}}{A^{np}\cdot B^{mq}}\frac{n!^{\mathrm{deg}(R_1)\epsilon}\cdot m!^{\mathrm{deg}(R_2)\epsilon}}{n!^p\cdot m!^q}
	$$
	Notice that $p+d_i=\mathrm{deg}(Q_i)$. Now for any $p,q\geq 1$, one can choose an $\epsilon>0$ such that $p-d_1\cdot \epsilon>9/10$ and $q-d_2\cdot \epsilon>9/10$. With such an $\epsilon$, we obtain 
	$$
	\frac{N(F(n,m))^{1+\epsilon}}{F(n,m)}\leq C''\frac{(4^{1+\epsilon})^n\cdot (4^{1+\epsilon})^m}{n!^{9/10}\cdot m!^{9/10}}\frac{1}{A^{(9/10)\cdot n}\cdot B^{(9/10)\cdot m}}
	$$
	Now use the Stirling approximation for $n!$ and $m!$ respectively, to conclude that 
	$$
	\frac{N(F(n,m))^{1+\epsilon}}{F(n,m)}\rightarrow 0
	$$
	as $(n,m)\rightarrow \infty$. This completes the proof.
	\noindent\\
	We now prove Proposition 1.7. We use the notation as in the proof of Theorem 1.1. We assume $n_r\geq \cdots \geq n_1$ and choose  $n_r>8\mathrm{max}\{A_1,...,A_r, f_{ij}\}$. According to Bertrand' postulate, there is a prime $q\in (n_r/2,n_r)$. For such a prime, we have
	$$
\nu_q((A_i^{n_i}n_i!)^{l_i}\cdot R_i(A_i^{n_i}n_i!))\leq l_i, 
$$ 
since $q$ divides neither $A_i$ nor $f_{i0}$. This implies 
 $$
 \nu_q(\prod_{i=1}^r(A_i^{n_i}n_i!)^{l_i}\cdot \prod_{i=1}^rR_i(A_i^{n_i}n_i!))\leq \sum_{i=1}^rl_i<d.
 $$
On the other hand we have 
$$
\nu_q(x^d)\geq d. 
$$
Since $d>l_1+\cdots +l_r$, it follows there are no solutions for $n_r>8\mathrm{max}\{A_1,...,A_r, f_{ij}\}$. This completes the proof.

	\vspace{0.3cm}
	\noindent
	{\tiny HOCHSCHULE FRESENIUS UNIVERSITY OF APPLIED SCIENCES 40476 D\"USSELDORF, GERMANY.}\\
	E-mail adress: sasa.novakovic@hs-fresenius.de\\
	
\end{document}